\documentclass{amsart}
\title{Digit patterns and Coleman power series}
\author{Greg W.\ Anderson}
\address{University of Minnesota, Mpls., MN 55455}
\email{gwanders@umn.edu}
\date{October 17, 2005}
\thanks{MSC-classification: 11S31.\\ \indent Keywords: Coleman power series, 
characteristic $p$, formal additive group, digit expansions}
\newcommand{\bracket}[1]{{\langle #1\rangle}}
\newcommand{\AAA}{{\mathcal{A}}}
\newcommand{\NN}{{\mathbb{N}}}
\newcommand{\MMM}{{\mathcal{M}}}
\newcommand{\NNN}{{\mathcal{N}}} 
\newcommand{\one}{{\mathbf{1}}}
\DeclareMathOperator{\ord}{{\mathrm{ord}}}
\newcommand{\iso}{{\stackrel{\sim}{\rightarrow}}}
\newcommand{\FF}{{\mathbb{F}}}
\newcommand{\ZZ}{{\mathbb{Z}}}
\newcommand{\QQ}{{\mathbb{Q}}}
\newcommand{\DDD}{{\mathbf{D}}}
\newtheorem{Proposition}[subsection]{Proposition}
\newtheorem{Theorem}[subsection]{Theorem}
\newtheorem{Corollary}[subsection]{Corollary}
\newtheorem{Lemma}[subsection]{Lemma}
\newcommand{\OO}{{\mathcal{O}}}
\newcommand{\FFF}{{\mathcal{F}}}
\begin{document}
\begin{abstract}
Our main result is elementary and concerns the relationship between the multiplicative groups of the coordinate and endomorphism rings of the formal additive group
over a field of characteristic $p>0$. The proof involves the combinatorics of base $p$ representations of positive integers in a striking way.
We apply the main result to construct a canonical quotient of the module of Coleman power series over the Iwasawa algebra 
when the base local field is of characteristic $p$. This gives information in a situation 
which apparently has never previously been investigated.
\end{abstract}
\maketitle
\section{Introduction}
\subsection{Overview and motivation}
Our main result (Theorem~\ref{Theorem:MainResult} below) concerns the relationship between the multiplicative groups of the coordinate and endomorphism rings of the formal additive group
over a field of characteristic $p>0$. Our result is elementary and does not require a great deal of apparatus for its statement. The proof of the main result involves the combinatorics of base $p$ representations of positive integers
in a striking way. 
We apply our main result (see Corollary~\ref{Corollary:ColemanApp} below)  to construct a canonical quotient of the module of Coleman power series over the Iwasawa algebra 
when the base local field is of characteristic $p$. By {\it Coleman power series} we mean the  telescoping power series
introduced and studied in Coleman's classical paper \cite{Coleman}. 

Apart from Coleman's \cite{ColemanLocModCirc} complete results in the important special case of the formal multiplicative group over $\ZZ_p$,  little is known about the structure of the module of Coleman  power series over the Iwasawa algebra, and, so far as we can tell, 
the characteristic $p$ situation has never previously been investigated. We undertook this research in an attempt to fill  the gap in characteristic $p$. 
Our results are far from being as complete as Coleman's, but they are surprising on account of their ``digital'' aspect, and they  raise further questions  worth investigating.

 \subsection{Formulation of the main result}   The notation introduced under this heading is in force throughout the paper.  
 
\subsubsection{Rings and groups of power series} Fix a prime number $p$ 
and a  field $K$ of characteristic $p$.
Let $q$ be a power of $p$. Consider: the (commutative) power series ring 
$$K[[X]]=\left\{\left.\sum_{i=0}^\infty a_i X^i\right| a_i\in K\right\};$$
the  (in general noncommutative) ring 
$$R_{q,K}=\left\{\left.\sum_{i=0}^\infty a_i X^{q^i}\right| a_i\in K\right\},$$
in which by definition multiplication is power series composition;
and the subgroup
$$\Gamma_{q,K}=\left\{\left.X+\sum_{i=1}^\infty a_i X^{q^i}\right| a_i\in K\right\}\subset R_{q,K}^\times,$$
where in general $A^\times$ denotes the group of units of a ring $A$ with unit.
Note that $K[[X]]^\times$ is a right $\Gamma_{q,K}$-module
via composition of power series. 

\subsubsection{Logarithmic differentiation}
Given $F=F(X)\in K[[X]]^\times$,
put 
$$\DDD[F](X)=XF'(X)/F(X)\in XK[[X]].$$
Note that
\begin{equation}\label{equation:Homogeneity}
\DDD[F(\alpha X)]=\DDD[F](\alpha X)
\end{equation}
for all $\alpha\in K^\times$. Note that the sequence
\begin{equation}\label{equation:Factoid}
1\rightarrow K[[X^p]]^\times\subset
K[[X]]^\times\xrightarrow{\DDD}
\left\{\left.\sum_{i=1}^\infty a_i X^i\in XK[[X]]\right|
a_{pi}=a_i^p\;\mbox{for all}\;i\in \NN \right\}\rightarrow 0
\end{equation} 
is exact,
where $\NN$ denotes the set of positive integers. 

\subsubsection{$q$-critical integers}
Given $c\in \NN$, let
$$
O_q(c)=\{n\in \NN\vert (n,p)=1\;\mbox{and}
\;n\equiv p^ic\bmod{q-1}\;\mbox{for some $i\in \NN\cup\{0\}$}\}.
$$
Given $n\in \NN$, let 
$\ord_p n$ denote the exact order with which $p$ divides $n$. 
We define
$$
C^0_q=\left\{c\in \NN\cap(0,q)\left|
(c,p)=1\;\mbox{and}\;\frac{c+1}{p^{\ord_p(c+1)}}=\min_{n\in O_q(c)\cap(0,q)}
\frac{n+1}{p^{\ord_p(n+1)}}\right.\right\},
$$
and we call elements of this set {\em $q$-critical integers}.
In the simplest case $p=q$ one has $C^0_p=\{1,\dots,p-1\}$, but in general
the set $C^0_q$ is somewhat complicated.
 Put
$$
C_q=\bigcup_{c\in C_q^0}\{q^i(c+1)-1\vert i\in \NN\cup\{0\}\},
$$
noting that the union is disjoint, since the sets in the union are contained in 
different congruence classes modulo $q-1$.
See below for informal ``digital'' descriptions of the sets $C_q^0$ and $C_q$.

\subsubsection{The homomorphism $\psi_q$}
We define a  homomorphism 
$$\psi_q:XK[[X]]\rightarrow X^2K[[X]]$$ as follows:
given  $F=F(X)=\sum_{i=1}^\infty a_iX^i\in XK[[X]]$,
put 
$$
\psi_{q}[F]=X\cdot \sum_{k\in C_q}a_kX^{k}.
$$
Note that the composite map
$$\psi_q\circ \DDD:K[[X]]^\times\rightarrow 
\left\{\left.\sum_{k\in C_q}a_k X^{k+1}\right| a_k\in K\right\}$$
is surjective by exactness of sequence (\ref{equation:Factoid}).
Further,  since the set $\{k+1\vert k\in C_q\}$ is stable under multiplication by $q$,
the target of $\psi_q\circ \DDD$ comes equipped with the structure of left $R_{q,K}$-module. More precisely, the target of $\psi_q\circ \DDD$ is a free left $R_{q,K}$-module for which the set $\{X^{k+1}\vert k\in C_q^0\}$ is a basis. 

The following is the main result  of the paper. 

\begin{Theorem}\label{Theorem:MainResult}
The formula
\begin{equation}\label{equation:MainResult}
\psi_{q}[\DDD[F\circ \gamma]]=\gamma^{-1}\circ \psi_{q}[\DDD[F]]
\end{equation}
holds for all $\gamma\in \Gamma_{q,K}$ and $F\in K[[X]]^\times$.
\end{Theorem}
\noindent 
 In \S\S\ref{section:Reduction}--\ref{section:DigitMadness2} we give the proof of the theorem. More precisely, 
we first explain in \S\ref{section:Reduction} how to reduce the proof of the theorem to a couple of essentially combinatorial assertions
(Theorems~\ref{Theorem:DigitMadness}
and \ref{Theorem:DigitMadness2}),
and then we prove the latter in
\S\ref{section:DigitMadness} and \S\ref{section:DigitMadness2}, respectively.  
In \S\ref{section:ColemanApp} we  make the application (Corollary~\ref{Corollary:ColemanApp}) of Theorem~\ref{Theorem:MainResult} to  Coleman power series. The application does not require any of the apparatus of the proof of Theorem~\ref{Theorem:MainResult}.

\subsection{Informal discussion} 
\subsubsection{``Digital'' description of $C_q^0$}
The definition of $C_q^0$ can readily be understood in terms
of simple operations on digit strings.
 For example, to verify that $39$ is $1024$-critical, begin by writing out the base $2$ representation of $39$ thus:
 $$39=100111_2$$
Then put enough place-holding $0$'s on the left so as to 
represent $39$ by a digit string of length $\ord_2 1024=10$:
$$39=0000100111_2$$
Then calculate as follows:
$$\begin{array}{cclr}
\mbox{permute cyclically}
&0000100111_2&\xrightarrow{\mbox{\tiny strike trailing $1$'s and leading $0$'s}}&100_2\\
&0001001110_2&\mbox{ignore: terminates with a $0$}\\
&0010011100_2&\mbox{ignore: terminates with a $0$}\\
\downarrow&0100111000_2&\mbox{ignore: terminates with a $0$}\\
&1001110000_2&\mbox{ignore: terminates with a $0$}\\
&0011100001_2&\xrightarrow{\mbox{\tiny strike trailing $1$'s and leading $0$'s}}&1110000_2\\
&0111000010_2&\mbox{ignore: terminates with a $0$}\\
&1110000100_2&\mbox{ignore: terminates with a $0$}\\
\downarrow
&1100001001_2&\xrightarrow{\mbox{\tiny strike trailing $1$'s and leading $0$'s}}&110000100_2\\
&1000010011_2&\xrightarrow{\mbox{\tiny strike trailing $1$'s and leading $0$'s}}&10000100_2\\
\end{array}$$
Finally, conclude that $39$ is $1024$-critical 
because the first entry of the last column is the smallest in that column. 
This numerical example conveys some of the flavor of the combinatorial considerations 
coming up in the proof of Theorem~\ref{Theorem:MainResult}. 
\subsubsection{``Digital'' description of $C_q$}
Given $c\in C_q^0$, let $[c_1,\cdots,c_m]_p$
be a string of digits representing $c$ in base $p$.
(The digit string notation is  defined below in
\S\ref{section:Reduction}.)
Then each digit string of the form
$$[c_1,\dots,c_m,\underbrace{p-1,\dots,p-1}_{ n\ord_p q}]_p$$
represents an element of $C_q$. Moreover, each element of $C_q$ arises this way, for unique $c\in C_q^0$ and $n\in \NN\cup\{0\}$.

\subsubsection{Miscellaneous remarks}
$\;$

(i) The set $C_q$ is a subset of the set of {\em magic numbers} (relative to the choice of $q$) as defined and studied in \cite[\S8.22, p.\ 309]{Goss}.   For the moment we do not understand this connection
 on any basis other than ``numerology'', but we suspect that  it runs much deeper.

(ii) A well-ordering of the set of positive integers distinct from the usual one,
which we call the {\em $p$-digital well-ordering},  plays a key role
in the proof of Theorem~\ref{Theorem:MainResult}, via Theorems~\ref{Theorem:DigitMadness}
and \ref{Theorem:DigitMadness2} below.
In particular, Theorem~\ref{Theorem:DigitMadness2}, via Proposition~\ref{Proposition:DigitMadness2},
characterizes the sets $C_q^0$ and $C_q$ in terms of the $p$-digital well-ordering
and congruences modulo $q-1$.

(iii)  The results of this paper were discovered by extensive
 computer experimentation with base $p$ expansions
 and binomial coefficients modulo $p$.  No doubt refinements of our results can be discovered 
 by continuing such experiments. 
 
 (iv) It is an open problem to find a minimal set of generators
 for $K[[X]]^\times$ as a topological right $\Gamma_{q,K}$-module, the topologies here being the $X$-adically induced ones.  It seems very likely that the module is always infinitely generated,
even when $K$ is a finite field.
Computer experimentation (based on the method of proof of Proposition~\ref{Proposition:Reduction} below) with the simplest case of the problem (in which $K$ is the two-element field and $p=q=2$) has revealed some interesting patterns. But still we are  unable to hazard any detailed guess about the solution.

\section{Application to Coleman power series}
\label{section:ColemanApp}
We assume that the reader already knows about Lubin-Tate formal groups and Coleman power series, and is familiar with their applications. We refer the less well-versed reader to \cite{LubinTate}, \cite{Coleman}, \cite{ColemanLocModCirc} and \cite{ColemanAI} to get started up the mountain of literature. 

\subsection{Background} We  review \cite{LubinTate}, \cite{Coleman} and \cite{ColemanLocModCirc} just far enough to fix a consistent system of notation and to frame precisely the general structure problem motivating our work. \subsubsection{The setting}
Let $k$ be a nonarchimedean local field with maximal
compact subring $\OO$ and uniformizer $\pi$.
Let $q$ and $p$ be the cardinality and characteristic, respectively,
of the residue field $\OO/\pi$. Let $\bar{k}$ be an algebraic closure of $k$.
Let $H$ be a complete unramified extension of $k$ in the completion of $\bar{k}$,
let $\varphi$ be the arithmetic Frobenius automorphism
of $H/k$, and let $\OO_H$ be the ring of integers of $H$. Let the action of $\varphi$ be
extended coefficient-by-coefficient to the power series
ring $\OO_H[[X]]$. 
\subsubsection{Lubin-Tate formal groups}
We say that formal power series with coefficients in $\OO$ are {\em congruent modulo $\pi$}
if they are so coefficient-by-coefficient,
and we say they are {\em congruent modulo degree $2$}
if the constant and linear terms agree. 
Let $\FFF_\pi$ be the set of one-variable power series $f=f(X)$ such that
$$f(X)\equiv \pi X\bmod{\deg 2},\;\;\;f(X)\equiv X^q\bmod{\pi}.$$
The simplest example of an element of $\FFF_\pi$ is 
$\pi X+X^q$.
The general example of an element of $\FFF_\pi$ is 
$\pi X+X^q+\pi X^2e(X)$,
where $e(X)\in \OO[[X]]$ is arbitrary.
Given $f\in \FFF_\pi$, there exists unique $F_f=F_f(X,Y)\in \OO[[X,Y]]$
such that
$$F_f(X,Y)\equiv X+Y\bmod{\deg 2},\;\;\;f(F_f(X,Y))=F_f(f(X),f(Y)).$$
The power series $F_f(X,Y)$ is a {\em commutative formal group law}.
Given $a\in \OO$ and $f,g\in \FFF_\pi$, there exists unique
$[a]_{f,g}=[a]_{f,g}(X)\in \OO[[X]]$ such that 
$$[a]_{f,g}(X)\equiv aX\bmod{\deg 2},\;\;f([a]_{f,g}(X))=[a]_{f,g}(g(X)).$$
We write $[a]_f=[a]_f(X)=[a]_{f,f}(X)$ to abbreviate notation.
The family \linebreak $\{[a]_f(X)\}_{a\in \OO}$ is a system of {\em formal complex multiplications}
for the formal group law $F_f(X,Y)$. For each fixed $f\in \FFF_\pi$, the package 
$$(\FFF_f(X,Y),\{[a]_f(X)\}_{a\in \OO})$$
is a {\em Lubin-Tate formal group}.
The formal properties of the ``big package''
$$\left(\{F_f(X,Y)\}_{f\in \FFF_\pi},
\{[a]_{f,g}(X)\}_{\begin{subarray}{c}
a\in \OO\\
f,g\in \FFF_\pi
\end{subarray}}\right)$$
are detailed in \cite[Thm.\ 1, p.\ 382]{LubinTate}. In particular,
one has
\begin{equation}\label{equation:LTformal}
[\pi]_f(X)=f(X),\;\;\;[1]_f(X)=X,\;\;\;
[a]_{f,g}\circ[b]_{g,h}=[ab]_{f,h}
\end{equation}
for all $a,b\in \OO$ and $f,g,h\in \FFF_\pi$.
We remark also that
\begin{equation}\label{equation:LTformalbis}
[\omega]_{\pi X+X^q}(X)=\omega X
\end{equation}
for all roots $\omega$ of unity in $k$ of order prime to $p$.

\subsubsection{Coleman power series}
 By Coleman's theory \cite{Coleman} there exists for each $f\in \FFF_\pi$ a unique
group homomorphism 
$$\NNN_f:\OO_H[[X]]^\times\rightarrow \OO_H[[X]]^\times$$
such that 
$$
\NNN_f[h](f(X))=\prod_{\begin{subarray}{c}
\lambda\in \bar{k}\\
f(\lambda)=0
\end{subarray}}h(F_f(X,\lambda))
$$
for all $h\in \OO_H[[X]]^\times$.
 Let 
$$\MMM_f=\{h\in \OO_H[[X]]^\times\vert \NNN_f[h]=\varphi h\}.$$
We refer to elements of $\MMM_f$ as {\em Coleman power series}. 

\subsubsection{Natural operations on Coleman power series}
The group
$\MMM_f$ comes equipped with the structure of right $\OO^\times$-module
by the rule 
\begin{equation}\label{equation:OOModuleRule}
((h,a)\mapsto h\circ [a]_f):\MMM_f\times \OO^\times\rightarrow \MMM_f,
\end{equation}
and we  have at our disposal a canonical isomorphism
\begin{equation}\label{equation:Cformal}
(h\mapsto h\circ [1]_{g,f}):\MMM_g\rightarrow \MMM_f
\end{equation}
of right $\OO^\times$-modules
for all $f,g\in \FFF_\pi$, as one verifies by applying the formal properties  (\ref{equation:LTformal}) of the big Lubin-Tate package in a straightforward way. We also have at our disposal a canonical group isomorphism
\begin{equation}\label{equation:CBijection}
(h\mapsto h\bmod{\pi}):\MMM_f\rightarrow (\OO_H/\pi)[[X]]^\times
\end{equation}
as one verifies by applying \cite[Lemma 13, p.\ 103]{Coleman}
in a straightforward way. 
The {\em Iwasawa algebra} (completed group ring)
$\ZZ_p[[\OO^\times]]$ associated to $k$ acts naturally on the slightly modified version 
$$\MMM_f^0=\{h\in \OO_H[[X]]^\times \vert h\in \MMM_f,\;h(0)\equiv 1\bmod{\pi}\}$$
of $\MMM_f$.

\subsubsection{The structure problem}
Little seems to be known in general about the structure of the $\OO^\times$-module $\MMM_f$. To determine this structure is a fundamental problem in local class field theory, and the problem remains open. Essentially everything we do know about the problem is due to Coleman.
Namely, in the special case
 $$k=\QQ_p=H,\;\;\;\pi=p,\;\;\;f(X)=(1+X)^p-1\in \FFF_\pi,$$
Coleman showed \cite{ColemanLocModCirc}
that $\MMM^0_f$ is ``nearly'' a free  $\ZZ_p[[\OO^\times]]$-module of rank $1$,
and in the process recovered Iwasawa's result on the structure of local units modulo circular units.  Moreover, Coleman's methods are strong enough to analyze $\MMM^0_f$ completely in the case of general $H$, even though this level of generality is not explicitly considered in \cite{ColemanLocModCirc}. So in the case of the formal multiplicative group over $\ZZ_p$
we have a complete and satisfying description of structure. Naturally one wishes for so
complete a description in the general case. We hope with the present work to contribute to the solution of the structure problem.

Here is the promised application of Theorem~\ref{Theorem:MainResult}, which makes the $\OO^\times$-module structure of a certain quotient of $\MMM_f$ explicit when $k$ is of characteristic $p$. 

\begin{Corollary}\label{Corollary:ColemanApp} Assume that $k$ is of characteristic $p$ and fix $f\in \FFF_\pi$.
Then there exists a surjective group homomorphism
$$\Psi_f:\MMM_f\rightarrow\left\{\left.\sum_{k\in C_q}a_kX^{k+1}\right| a_k\in \OO_H/\pi\right\}$$
such that
\begin{equation}\label{equation:ColemanApp}
\Psi_f[h\circ [\omega u]_f]\equiv
[(\omega u)^{-1}]_{\pi X+X^q}\circ \Psi_f[h]\circ[\omega]_{\pi X+X^q}\bmod{\pi}
\end{equation}
for all $h=h(X)\in \MMM_{f}$, $u\in 1+\pi\OO$
and roots of unity $\omega\in \OO^\times$.
 \end{Corollary}
\proof If we are able to construct $\Psi_{\pi X+X^q}$ with the desired properties,
then in the general case the map
$$\Psi_f=(h\mapsto \Psi_{\pi X+X^q}[h\circ [1]_{f,\pi X+X^q}])$$
has  the desired properties by (\ref{equation:LTformal}) and (\ref{equation:Cformal}).
We may therefore assume without loss of generality that 
$$f=\pi X+X^q,$$
 in which case 
 $F_f(X,Y)=X+Y$,
 i.~e., the formal group underlying the Lubin-Tate formal group attached to $f$ is additive.
 By (\ref{equation:LTformalbis}) and the definitions, given
 $a\in \OO$ and writing
 $$a=\sum_{i=0}^\infty \alpha_i \pi^i\;\;\;(\alpha_i^q=\alpha_i),$$
 in the unique possible way, one has
 $$[a]_f\equiv \sum_{i=0}^\infty \alpha_i X^{q^i}\bmod{\pi},$$
 and hence 
the map
 $a\mapsto [a]_{f}\bmod \pi$ gives rise to an isomorphism
$$\theta:\OO\iso R_{q,\OO/\pi}$$
of rings. Let
$$\rho:\MMM_f\rightarrow (\OO_H/\pi)[[X]]^\times$$
be the isomorphism (\ref{equation:CBijection}). We claim that
$$\Psi_f=\psi_q\circ \DDD\circ \rho$$
has all the desired properties. In any case, since $\rho$ is an isomorphism and $\psi_q\circ \DDD$ by (\ref{equation:Factoid}) is surjective, $\Psi_f$  is surjective, too.
To verify (\ref{equation:ColemanApp}), we calculate as follows:
$$\begin{array}{rcl}
\psi_q[\DDD[\rho(h\circ [\omega u]_f)]]&=&\psi_q[\DDD[\rho(h\circ [\omega]_f\circ[ u]_f)]]\\
&=&
\psi_q[\DDD[\rho(h)\circ \theta(\omega)\circ \theta(u)]]\\
&=&\theta(u^{-1})\circ \psi_q[\DDD[\rho(h)\circ \theta(\omega)]]\\
&=&\theta(u^{-1})\circ \psi_q[\DDD[\rho(h)]\circ \theta(\omega)]\\
&=&\theta(u^{-1})\circ\theta(\omega^{-1})\circ
\psi_q[\DDD[\rho(h)]]\circ \theta(\omega)\\
&=&\theta((u\omega)^{-1})\circ
\psi_q[\DDD[\rho(h)]]\circ \theta(\omega)
\end{array}
$$
The third and fourth steps are justified by
(\ref{equation:MainResult}) and 
(\ref{equation:Homogeneity}), respectively.
 The remaining steps are clear. The claim is proved, and with it the corollary. 
 \qed

\section{Reduction of the proof}\label{section:Reduction} We put Coleman power series behind us for the rest of the paper. We return to the elementary point of view taken in the introduction. In this section
we explain how to reduce the proof of Theorem~\ref{Theorem:MainResult} to a couple of combinatorial assertions.
\subsection{Digital apparatus}

\subsubsection{Base $p$ expansions}
Given an additive decomposition
$$n=\sum_{i=1}^s n_ip^{s-i}\;\;\;(n_i\in \ZZ\cap [0,p),\;n\in \NN),$$
we write
$$n=[n_1,\dots,n_s]_p,$$
we call the latter a {\em base $p$ expansion} of $n$ and 
we call the coefficients $n_i$ {\em digits}.
Note that we allow base $p$ expansions to have leading $0$'s.
We say that a base $p$ expansion is
{\em minimal} if the first digit is positive.
For convenience, we set the empty base $p$ expansion $[]_p$ equal to $0$ and declare it  to be minimal.
 We always read base $p$ expansions left-to-right, as though they were words spelled in the alphabet $\{0,\dots,p-1\}$.  In this notation the well-known theorem of Lucas takes the form
$$\left(\begin{array}{c}
\;[a_1,\dots,a_n]_p\\
\;[b_1,\dots,b_n]_p
\end{array}\right)\equiv \left(\begin{array}{c}
a_1\\
b_1
\end{array}\right)
\cdots \left(\begin{array}{c}
a_n\\
b_n
\end{array}\right)\bmod{p}.$$
(For all $n\in \NN\cup\{0\}$ and $k\in \ZZ$  we set
$\left(\begin{subarray}{c}
n\\
k\end{subarray}\right)=
\frac{n!}{k!(n-k)!}$ if $0\leq k\leq n$ and $\left(\begin{subarray}{c}
n\\
k\end{subarray}\right)=0$ otherwise.)
The theorem of Lucas implies  that for all 
integers $k,\ell,m\geq 0$ such that $m=k+\ell$, the binomial
coefficient $\left(\begin{subarray}{c}
m\\
k
\end{subarray}\right)$ does not vanish modulo $p$ if and only if the addition of $k$ and $\ell$ in base $p$ requires no ``carrying''.
\subsubsection{The $p$-core function $\kappa_p$}
Given $n\in \NN$, we define
$$
\kappa_p(n)=(n/p^{\ord_p n}+1)/p^{\ord_p(n/p^{\ord_p n}+1)}-1.
$$
We call $\kappa_p(n)$  the {\em $p$-core}
 of $n$.
For example, $\kappa_p(n)=0$
iff  $n=p^{k-1}(p^\ell-1)$
for some $k,\ell\in \NN$.   The meaning of the $p$-core function
is easiest to grasp in terms of minimal base $p$ expansions. One calculates $\kappa_p(n)$ by 
discarding trailing $0$'s and then discarding trailing $(p-1)$'s.
For example, to calculate the $3$-core
of
$963=[1,0,2,2,2,0,0]_3$,
first discard trailing $0$'s to get
$[1,0,2,2,2]_3=107$,
and then discard  trailing $2$'s to get
$\kappa_3(963)=[1,0]_3=3$. \subsubsection{The $p$-defect function $\delta_p$}
For each $n\in\NN$, let $\delta_p(n)$ be the length of the minimal base $p$
representation of $\kappa_p(n)$. We call $\delta_p(n)$ the {\em $p$-defect}
of $n$. For example, since as noted above $\kappa_3(963)=[1,0]_3$,
one has $\delta_3(963)=2$.

\subsubsection{The $p$-digital well-ordering}
We equip the set of positive integers with a well-ordering 
 $\leq_p$ by declaring $m\leq_{p}n$
 if
 $$ \kappa_p(m)<\kappa_p(n)$$
 or
 $$\kappa_p(m)=\kappa_p(n)\;\mbox{and}\;m/p^{\ord_p m}<n/p^{\ord_pn}$$
 or
 $$\kappa_p(m)=\kappa_p(n)\;\mbox{and}\;m/p^{\ord_p m}=n/p^{\ord_p n}\;\mbox{and} \;m\leq n.$$
 In other words,  to verify $m\leq _p n$, first compare $p$-cores of $m$ and $n$, 
then in case of a tie compare numbers of $(p-1)$'s trailing the $p$-core, and in case of  another tie compare numbers of trailing $0$'s. We call $\leq_p$ the {\em $p$-digital well-ordering}.
 In the obvious way we derive order relations $<_p$, $\geq_{p}$ and $>_{p}$
 from $\leq_{p}$. We remark that
 $$\delta_p(m)<\delta_p(n)\Rightarrow m<_p n,\;\;\;
 m\leq_p n\Rightarrow \delta_p(m)\leq \delta_p(n);$$
in other words, the function $\delta_p$ gives a reasonable if rough approximation
 to the $p$-digital well-ordering.
\subsubsection{The function $\mu_q$}
Given $c\in \NN$, let $\mu_q(c)$ be the unique element of the set 
$$\{n\in \NN\vert n\equiv p^i c\bmod{q-1}\;\mbox{for some}\;i\in \NN\cup\{0\}\}
$$
minimal with respect to the $p$-digital well-ordering.
Note that $\mu_q(c)$ cannot be divisible by $p$.
Consequently $\mu_q(c)$ may also be characterized as the unique element of the set $O_q(c)$ minimal with respect to the $p$-digital well-ordering. 
\subsubsection{$p$-admissibility}
We say that a quadruple $(j,k,\ell,m)\in \NN^4$
is {\em $p$-admissible} if
$$(m,p)=1,\;\;\;m=k+j(p^\ell-1),\;\;\;\left(\begin{array}{c}
k-1\\
j
\end{array}\right)\not\equiv 0\bmod{p}.$$
This is the key technical definition of the paper.
Let $\AAA_p$ denote the set of $p$-admissible quadruples.

\begin{Theorem}\label{Theorem:DigitMadness}
For all $(j,k,\ell,m)\in \AAA_p$, one has (i)
$k<_p m$, and moreover, (ii)
if \linebreak $\kappa_p(k)=\kappa_p(m)$, then  
 $j=(p^{\ord_p k}-1)/(p^\ell-1)$.
 \end{Theorem}
\noindent We will prove this result in \S\ref{section:DigitMadness}. Note that the conclusion of part (ii) of the theorem
implies $\ord_pk>0$ and $\ell\vert \ord_p k$.

\begin{Theorem}\label{Theorem:DigitMadness2} 
One has
\begin{equation}\label{equation:DigitMadness2}
\max_{c\in \NN}\mu_q(c)<q,
\end{equation}
\begin{equation}\label{equation:DigitMadness2bis}
\begin{array}{cl}
&\displaystyle\{(\mu_q(c)+1)q^{i}-1\;\vert\; i\in \NN\cup\{0\},\;c\in \NN\}\\\\
=&\displaystyle\left\{c\in \NN\left|(c,p)=1,\;
\kappa_p(c)=\min_{n\in O_q(c)}\kappa_p(n)\right.\right\}.
\end{array}
\end{equation}
\end{Theorem}
\noindent
 We will prove this result in \S\ref{section:DigitMadness2}.  We have phrased the result in a way emphasizing the $p$-digital well-ordering.
 But perhaps it is not clear what the theorem means in the context of Theorem~\ref{Theorem:MainResult}.
The next result provides an explanation.

 \begin{Proposition}\label{Proposition:DigitMadness2}
 Theorem~\ref{Theorem:DigitMadness2} granted,
 one has
 \begin{equation}\label{equation:DigitMadness2quad}
 C_q^0=\{\mu_q(c)\vert c\in \NN\},
 \end{equation}
 \begin{equation}\label{equation:DigitMadness2ter}
 C_q=\left\{c\in \NN\left|(c,p)=1,\;
\kappa_p(c)=\kappa_p(\mu_q(c))\right.\right\}.
 \end{equation}
 \end{Proposition}
 \proof
 The definition of $C^0_q$ can be rewritten $$C^0_q=\left\{c\in \NN\cap(0,q)\left|
 (c,p)=1,\;\kappa_p(c)=\min_{n\in O_q(c)\cap(0,q)}\kappa_p(n)\right.\right\}.$$
Therefore relation (\ref{equation:DigitMadness2}) implies containment $\supset$ in (\ref{equation:DigitMadness2quad}) and moreover, supposing failure of equality in (\ref{equation:DigitMadness2quad}), there exist
$c,c'\in C_q^0$
such that 
$$c=\mu_q(c)\neq c',\;\;\;\kappa_p(c)=\kappa_p(c').$$
But $c'=q^i(c+1)-1$ for some $i\in \NN$ by (\ref{equation:DigitMadness2bis}), hence $c'\geq q$,
and hence $c'\not\in C_q^0$. This contradiction establishes equality in (\ref{equation:DigitMadness2quad})
 and in turn containment $\subset$ in (\ref{equation:DigitMadness2ter}).
Finally, (\ref{equation:DigitMadness2bis})
 and (\ref{equation:DigitMadness2quad}) imply equality in (\ref{equation:DigitMadness2ter}). \qed

The following is the promised reduction of the proof of Theorem~\ref{Theorem:MainResult}.

\begin{Proposition}\label{Proposition:Reduction}
If Theorems~\ref{Theorem:DigitMadness} 
and \ref{Theorem:DigitMadness2} hold, then 
Theorem~\ref{Theorem:MainResult} holds, too.
\end{Proposition}
\noindent Before turning to the proof, we pause to discuss the groups in play.
\subsection{Generators for $K[[X]]^\times$,
$\DDD[K[[X]]^\times]$ and $\Gamma_{q,K}$}
\label{subsection:Convenient}
Equip $K[[X]]^\times$ with the topology for which the family 
$\{1+X^nK[[X]]\vert n\in \NN\}$ is a neighborhood base at the origin. Then the set
$$\{1+\alpha X^k\vert \alpha \in K^\times,\;k\in \NN\}\cup K^\times$$
generates $K[[X]]^\times$ as a topological group.
Let $\FF_p$ be the residue field of $\ZZ_p$.
Let 
$E_p=E_p(X)\in \FF_p[[X]]$
be the reduction modulo $p$ of the Artin-Hasse exponential
$$\exp\left(\sum_{i=0}^\infty
\frac{X^{p^i}}{p^i}\right)\in (\QQ\cap\ZZ_p)[[X]],$$
noting that
$$\DDD[E_p]=\sum_{i=0}^\infty X^{p^i}.$$
Since
$E_p(X)=1+X+O(X^2)$,
the set
$$\{E_p(\alpha X^k)\;\vert \;\alpha\in K^\times,\;k\in \NN,\;(k,p)=1\}\cup K[[X^p]]^\times$$
generates $K[[X]]^\times$ as a topological group.
For each $k\in \NN$ such that $(k,p)=1$ and $\alpha\in K^\times$,  put
$$W_{k,\alpha}=W_{k,\alpha}(X)=k^{-1}\DDD[E_p(\alpha X^k)]=\sum_{i=0}^\infty \alpha^{p^i}X^{kp^i}\in XK[[X]].
$$
Equip $\DDD[K[[X]]^\times]$ with the relative $X$-adic topology. The set
$$\{W_{k,\alpha}\vert k\in \NN,\;(k,p)=1,\;\alpha\in K^\times\}$$
generates $\DDD[K[[X]]^\times]$ as a topological group, cf.\ exact sequence (\ref{equation:Factoid}).
Equip $\Gamma_{q,K}$ with the relative $X$-adic topology.
Note that
$$
(X+\beta X^{q^\ell})^{-1}=\sum_{i=0}^\infty
(-1)^{i}\beta^{\frac{q^{\ell i}-1}{q^\ell-1}}X^{q^{\ell i}}\in \Gamma_{q,K}
$$
for all $\ell\in \NN$ and $\beta \in K^\times$. The inverse operation here is of course understood in the functional rather than multiplicative sense.
The set
$$\{X+\beta X^{q^\ell}\;\vert\;
\beta\in K^\times,\;\ell\in \NN\}$$
generates $\Gamma_{q,K}$ as a topological group.
\subsection{Proof of the proposition}
It is enough to verify (\ref{equation:MainResult}) with $F$ and $\gamma$ ranging over sets of generators for the topological groups $K[[X]]^\times$
and $\Gamma_{q,K}$, respectively. The generators mentioned in the preceding paragraph are the convenient ones.  So fix $\alpha,\beta\in K^\times$
and $k,\ell\in \NN$ such that $(k,p)=1$.   It will be enough to verify that
\begin{equation}\label{equation:Nuff}
\psi_{q}[M_{k,\alpha,\ell,\beta}]=
\left\{\begin{array}{rl}\displaystyle
\alpha X^{k+1}+\sum_{\ell\vert f\in \NN} (-1)^{f/\ell}\alpha^{q^f}\beta^{\frac{q^f-1}{q^\ell-1}}X^{q^f(k+1)}&\mbox{if $k\in C_q$,}\\\\
0&\mbox{otherwise,}
\end{array}\right.
\end{equation}
where  \begin{equation}\label{equation:MExpansion}
\begin{array}{cl}
&M_{k,\alpha,\ell,\beta}=M_{k,\alpha,\ell,\beta}(X)=k^{-1}\DDD[E_p(\alpha(X+\beta X^{q^\ell})^k)]\\\\
=&\displaystyle W_{k,\alpha}+\sum_{i=0}^\infty
\sum_{j=1}^\infty \left(\begin{array}{c}
p^ik-1\\ j
\end{array}\right)\alpha^{p^i}\beta^{j}X^{p^ik+j(q^\ell-1)}.\end{array}
\end{equation}
By Theorem~\ref{Theorem:DigitMadness},  many terms
on the right side of (\ref{equation:MExpansion}) vanish, and more  precisely,
we can rewrite (\ref{equation:MExpansion}) as follows:
\begin{equation}\label{equation:Nuff2}
\begin{array}{rcl}
M_{k,\alpha,\ell,\beta}&\equiv&\displaystyle\alpha X^k+
\sum_{\begin{subarray}{c}
m\in O_q(k)\\
m>_pk
\end{subarray}}
\left(\sum_{\begin{subarray}{c}
i\in \NN\cup\{0\}, j\in\NN\\
(j,p^ik,\ord_p q^\ell,m)\in \AAA_p
\end{subarray}}
\left(\begin{array}{c}
p^ik-1\\
j
\end{array}\right)\alpha^{p^i}\beta^j\right)X^m\\\\
&&\hskip 5cm\bmod{X^pK[[X^p]]}.
\end{array}
\end{equation}
By Theorem~\ref{Theorem:DigitMadness2} as recast in the form of Proposition~\ref{Proposition:DigitMadness2}, along with formula (\ref{equation:Nuff2}) and the definitions, both sides of (\ref{equation:Nuff})
vanish unless $k\in C_q$. So now fix $c\in C_q^0$ and $g\in \NN\cup\{0\}$
and put 
$$k=(c+1)q^g-1\in C_q$$ for the rest of the proof of the proposition.
Also fix $f\in \NN\cup\{0\}$ and put 
$$m=q^f(k+1)-1=(c+1)q^{f+g}-1\in C_q$$ for the rest of the proof. 
It is enough to evaluate the coefficient
of $X^m$ in (\ref{equation:Nuff2}). By part (ii) of Theorem~\ref{Theorem:DigitMadness}, there is no term in the sum on the right side of (\ref{equation:Nuff2}) of degree $m$ unless $\ell\vert f$, in
which case there is exactly one term, namely
$$\left(\begin{array}{c}
q^fk-1\\
\frac{q^f-1}{q^\ell-1}
\end{array}\right)\alpha^{q^f}\beta^{\frac{q^f-1}{q^\ell-1}}X^m,$$
and by the theorem of Lucas, the binomial coefficient mod $p$ evaluates 
 to $(-1)^{f/\ell}$.
Therefore (\ref{equation:Nuff}) does indeed hold.
\qed
\subsection{Remarks}
$\;$

(i) By formula (\ref{equation:Nuff2}), the $p$-digital well-ordering actually gives rise to a  $\Gamma_{q,K}$-stable complete separated filtration of the quotient $K[[X]]^\times/K[[X^p]]^\times$ distinct from the $X$-adically induced one. Theorem~\ref{Theorem:MainResult} merely describes the structure of $K[[X]]^\times/K[[X^p]]^\times$ near the top of the ``$p$-digital filtration''.

(ii) Computer experimentation based on formula (\ref{equation:MExpansion}) was helpful in making the discoveries detailed in this paper. We believe that continuation of such experiments could lead to further progress, e.g., to the discovery of a minimal set of generators for $K[[X]]^\times$ as a topological right $\Gamma_{q,K}$-module.

\section{Proof of Theorem~\ref{Theorem:DigitMadness}}\label{section:DigitMadness}
 \begin{Lemma}\label{Lemma:DigitGames}
Fix $(j,k,\ell,m)\in \AAA_p$. Put
$$e=\ord_p(m+1),\;\;\;f=\ord_p k,\;\;\;g=\ord_p(k/p^f+1).$$
Then there exists a unique integer $r$ such that
\begin{equation}\label{equation:DigitGames0}
0\leq r\leq e+\ell-1,\;\; r\equiv 0\bmod{\ell},\;\; j\equiv \frac{p^r-1}{p^\ell-1}\bmod{p^{e}},
\end{equation}
and moreover 
\begin{equation}\label{equation:DigitGames1}
f+g\geq e,
\end{equation}
\begin{equation}\label{equation:DigitGames2}
\kappa_p(m)\geq \kappa_p(k).
\end{equation}
\end{Lemma}
\noindent This lemma is the key technical result of the paper.
 \subsection{Completion of the proof of the theorem, granting the lemma}
Fix $(j,k,\ell,m)\in \AAA_p$.  
Let $e,f,g,r$ be as defined in Lemma~\ref{Lemma:DigitGames}.
Since the number of digits in the minimal base $p$ expansion of $k$ cannot exceed the number of digits in the minimal base $p$ expansion of $m$, one has
\begin{equation}\label{equation:DigitGames3}
\delta_p(k)+f+g\leq \delta_p(m)+e.
\end{equation}
Combining (\ref{equation:DigitGames1})
and (\ref{equation:DigitGames3}), one has
\begin{equation}\label{equation:DigitGames4}
\delta_p(k)=\delta_p(m)\Rightarrow f+g=e.
\end{equation}
Now in general one has
$$m+1=(\kappa_p(m)+1)p^e,\;\;\;
k+p^f=(\kappa_p(k)+1)p^{f+g},$$
and hence 
$$
\kappa_p(k)=\kappa_p(m)\Rightarrow 
\left(j=\frac{p^{f}-1}{p^\ell-1}\;\mbox{and}\;e>g\right)
$$
via 
(\ref{equation:DigitGames4}).
Theorem~\ref{Theorem:DigitMadness} now follows via (\ref{equation:DigitGames2})
and the definition of the $p$-digital well-ordering.
\qed

\subsection{Proof of Lemma~\ref{Lemma:DigitGames}}
Since $e$ is the number of trailing $(p-1)$'s in the minimal base $p$ expansion of $m$, the lemma is trivial in the case $e=0$.
We therefore assume that $e>0$ for the rest of the proof.

  Let 
$$m=[m_1,\dots,m_t]_p\;\;\;(t>0,\;m_1>0,\;\;m_t>0)$$
be the minimal base $p$ expansion of $m$. For convenience, put
$$d=\delta_p(m)\geq 0,\;\;\;m_\nu=0\;\mbox{for $\nu<1$}.$$
Then
$$t=e+d,\;\;\;m_{d+1}=\cdots=m_{d+e}=p-1,\;\;\;m_{d}< p-1.$$
By hypothesis
$$\left(\begin{array}{c}
k-1\\
j
\end{array}\right)=\left(\begin{array}{c}
m-jp^\ell-1+j\\
m-jp^\ell-1
\end{array}\right)>0,$$
hence
$$m>jp^\ell,$$
and hence the number of digits in the minimal base $p$ of expansion of $jp^\ell$ does not exceed that of $m$. Accordingly,
$$t> \ell$$
and one has a base $p$ expansion
for $j$ of the form
$$j=[j_1,\dots,j_{t-\ell}]_p,$$
which perhaps is not minimal. 
For convenience, put 
$$j_\nu=0\;\mbox{for $\nu<1$ and also for $\nu>t-\ell$}.$$
This state of affairs is summarized by the ``snapshot''
$$m=[m_1,\dots,m_t]=[m_{1},\dots,m_{d},\underbrace{p-1,\dots,p-1}_e]_p,\;\;\kappa_p(m)=[m_{1},\dots,m_{d}]_p,$$
$$jp^\ell=[j_1,\dots,j_t]_p=[j_1,\dots,j_{t-\ell},\underbrace{0,\dots,0}_\ell]_p
,$$
which the reader should keep in mind as we proceed. 

We are ready now to prove the existence and uniqueness of $r$. One has
$$m-jp^\ell-1=k-1-j=[m_1',\dots,m_d',p-1-j_{d+1},\dots,
p-1-j_{t-1},p-2]_p,$$
where the digits $m'_1,\dots,m'_d$ are defined by the equation
\begin{equation}\label{equation:Swivel}
\kappa_p(m)-[j_{1},\dots,j_{d}]_p=[m'_1,\dots,m'_d]_p.
\end{equation}
By hypothesis and the theorem of Lucas, the addition of $k-1-j$ and $j$ in base $p$ requires no ``carrying'', and hence
\begin{equation}\label{equation:BigDigit}
k-1=
\left\{\begin{array}{ll}
\;[m_1'+j_{1-\ell},\dots,m_d'+j_{d-\ell},\\
\;p-1-j_{d+1}+j_{d+1-\ell},\dots,
p-1-j_{d+e-1}+j_{d+e-1-\ell},p-2+j_{d+e-\ell}]_p.
\end{array}\right.
\end{equation}
From the system of inequalities for the last $e+\ell$ digits of the base $p$ expansion of $jp^\ell$
implicit in (\ref{equation:BigDigit}), it follows that there exists $r_0\in \NN\cup\{0\}$ such that
\begin{equation}\label{equation:BigDigitBis}
jp^\ell=[j_{1-\ell},\dots,j_{d-\ell},
\overbrace{0,\dots,0,\underbrace{\underbrace{1,0,\dots,0}_{\ell},\dots,\underbrace{1,0,\dots,0}_{\ell}}_{\mbox{\tiny $r_0$ blocks}},0}^{e+\ell}]_p.
\end{equation}
Therefore $r=r_0\ell$ has the required properties (\ref{equation:DigitGames0}). Uniqueness of $r$ is clear.   For later use, note the relation
\begin{equation}\label{equation:HeadScratcher}
r\geq e\Leftrightarrow [j_{d-\ell+1},\dots,j_{d}]_p\neq 0\Rightarrow [j_{d-\ell+1},\dots,j_{d}]_p=p^{r-e},
\end{equation}
which is  easy to see from the point of view adopted here to prove (\ref{equation:DigitGames0}).

By (\ref{equation:DigitGames0}) one has
\begin{equation}\label{equation:DigitGames2.2}
k+p^r-(m+1)+j' p^e (p^\ell-1)=0\;\;\;\mbox{for some $j'\in \NN\cup\{0\}$,}
\end{equation}
and hence one has
\begin{equation}\label{equation:DigitGames2.5}
r\geq \min(f,e),\;\;\;f\geq \min(r,e).
\end{equation}
This proves   (\ref{equation:DigitGames1}),
since either one has $f\geq e$, in which case (\ref{equation:DigitGames1}) holds trivially,
or else $f<e$, in which case $r=f$ by (\ref{equation:DigitGames2.5}), and hence
(\ref{equation:DigitGames1}) holds by (\ref{equation:DigitGames2.2}).

Put
$$k-1=[k'_1,\dots,k'_{d+e}]_p,\;\;\;\;\one_{r\geq e}=\left\{\begin{array}{ll}
1&\mbox{if $r\geq e$,}\\
0&\mbox{if $r<e$.}
\end{array}\right.
$$
Comparing (\ref{equation:BigDigit}) and (\ref{equation:BigDigitBis}), we see that 
 the digits $k'_{d+1},\dots,k'_{d+e}$ are all $(p-1)$'s with at most one exception, 
and the exceptional digit if it exists is a $p-2$.
Further, one has
$$k'_{d+1}=\dots=k'_{d+e}=p-1\Leftrightarrow f\geq e\Leftrightarrow \one_{r\geq e}=1$$
by (\ref{equation:DigitGames2.5}).
Therefore one has
$$\kappa_p(k)\leq[k'_1,\dots,k'_d]+\one_{r\geq e}.$$
Finally, via (\ref{equation:Swivel}), (\ref{equation:BigDigit}) and (\ref{equation:HeadScratcher}), it follows that
$$
\begin{array}{rcl}
\kappa_p(k)&\leq &[m_1'+j_{1-\ell},\dots,m_d'+j_{d-\ell}]_p+\one_{r\geq e}\\
&=&\kappa_p(m)-[j_1,\dots,j_d]_p+[j_{1-\ell},\dots,j_{d-\ell}]_p+\one_{r\geq e}\\
&=&\kappa_p(m)-[j_{1-\ell},\dots,j_d]_p+[j_{1-\ell},\dots,j_{d-\ell}]_p+\one_{r\geq e}\\
&=&\kappa_p(m)-[j_{d-\ell+1},\dots,j_{d}]_p+\one_{r\geq e}
\\
&&-[j_{1-\ell},\dots,j_{d-\ell},\underbrace{0,\dots,0}_\ell]_p+[j_{1-\ell},\dots,j_{d-\ell}]_p\\
&=&\kappa_p(m)-\one_{r\geq e}(p^{r-e}-1)-(p^\ell-1)[j_{1-\ell},\dots,j_{d-\ell}]_p\\
&\leq &\kappa_p(m).
\end{array}
$$
Thus (\ref{equation:DigitGames2}) holds and the proof of the lemma is complete.
\qed

\section{Proof of Theorem~\ref{Theorem:DigitMadness2}}\label{section:DigitMadness2}

\subsection{Further digital apparatus}
Put $\lambda=\ord_p q$.
For each $c\in \NN$, let
$$\bracket{c}_q=\min\{n\in \NN\vert n\equiv c\bmod{q-1}\},\;\;\;\tau_p(c)=c/p^{\ord_pc}.
$$
Note that
$$0<\bracket{c}_q<q,\;\;\;
\bracket{c}_q=\bracket{c'}_q\Leftrightarrow c\equiv c'\bmod{q-1}$$
for all $c,c'\in \NN$.
Given $c\in \NN$, and writing $\bracket{c}_q=[c_1,\dots,c_\lambda]_p$,
note that
$$\{c_1,\dots,c_\lambda\}\neq \{0\},\;\;
\bracket{pc}_q=[c_2,\dots,c_\lambda,c_1]_p,$$
$$\langle c\rangle_q\geq \tau_p(\bracket{c}_q)=[c_1,\dots,c_{\max\{i\vert c_i\neq 0\}}]_p\geq \kappa_p(\langle c\rangle_q).$$

\begin{Lemma}\label{Lemma:Necklace2}
$\langle p^ic\rangle_q\leq
p^i-1\Rightarrow \tau_p(\langle c\rangle_q)\leq\langle p^ic\rangle_q$ for $c\in \NN$ and $i\in \NN\cap(0,\lambda)$.
\end{Lemma}
\begin{Lemma}\label{Lemma:Necklace3}
$\displaystyle 
\min_{i=0}^{\lambda-1}\tau_p(\langle p^ic+1\rangle_q)=1+\min_{i=0}^{\lambda-1}\kappa_p(\langle p^ic\rangle_q)=1+\kappa_p(\mu_q(c))$ for $c\in \NN$.
\end{Lemma}
\begin{Lemma}\label{Lemma:Necklace4}
 $i\not\equiv 0\bmod{\lambda}\Rightarrow p^i(\mu_q(c)+1)-1\not\in O_q(c)$  for $i,c\in \NN$.
 \end{Lemma}
\subsection{Completion of the proof of the theorem, granting the lemmas}
Relation (\ref{equation:DigitMadness2}) holds
by Lemma~\ref{Lemma:Necklace3}.
Relation (\ref{equation:DigitMadness2bis}) holds
by Lemma~\ref{Lemma:Necklace4}.
\qed

\subsection{Proof of Lemma~\ref{Lemma:Necklace2}} Write $\langle c\rangle_q=[c_1,\dots,c_\lambda]_p$.
By hypothesis
$$\langle p^ic\rangle_q=
[\underbrace{0,\dots,0}_{\lambda-i},c_1,\dots,c_i]_q,\;\;\;
c=[c_1,\dots,c_i,\underbrace{0,\dots,0}_{\lambda-i}]_p,
$$
and hence $\tau_p(c)\leq \langle p^ic\rangle_q$.
\qed

\subsection{Proof of Lemma~\ref{Lemma:Necklace3}}
Since
$$\mu_q(c)=(\kappa_p(\mu_q(c))+1)p^g-1\in O_q(c),$$
for some $g\in \NN\cup\{0\}$,
one has
$$\kappa_p(\mu_q(c))+1\geq 
\min_{i=0}^{\lambda-1}
\min_{j=0}^{\lambda-1}
\langle p^i(p^jc+1)\rangle_q.$$
One has
$$
\tau_p(\langle n+1\rangle_q)\geq 1+\kappa_p(\langle n\rangle_q)
$$
for all $n\in \NN$, as can be verified by a somewhat tedious case analysis
which we omit. Clearly, the inequalities $\geq$  hold in the statement we are trying to prove.
Therefore it will be enough to prove that
$$
\min_{i=0}^{\lambda-1}
\min_{j=0}^{\lambda-1}
\langle p^i(p^jc+1)\rangle_q\geq \min_{j=0}^{\lambda-1} \tau_p(\langle p^jc+1\rangle_q).
$$
Fix $i=1,\dots,\lambda-1$ and $j=0,\dots,\lambda-1$.
It will be enough just to prove that
\begin{equation}\label{equation:AlmostLastNuff}
\langle p^i(p^jc+1)\rangle_q
<\tau_p(\langle p^jc+1\rangle_q)
\Rightarrow \langle p^i(p^jc+1)\rangle_q
\geq \tau_p(\langle p^{i+j}c+1\rangle_q).
\end{equation}
But by the preceding lemma, under the hypothesis of (\ref{equation:AlmostLastNuff}), one has 
$$p^i-1<\langle p^i(p^jc+1)\rangle_q$$
and hence
$$\langle p^i(p^jc+1)\rangle_q=
\langle p^{i+j}c+1\rangle_q+p^i-1\geq \tau_p(
\langle p^{i+j}c+1\rangle_q).$$
Thus (\ref{equation:AlmostLastNuff}) is proved, and with
it the lemma.
 \qed

\subsection{Proof of Lemma~\ref{Lemma:Necklace4}}
We may assume without loss of generality that \linebreak $0<i<\lambda$ and $c=\mu_q(c)$.
By the preceding lemma $c<q$.
Write
$c=[c_1,\dots,c_\lambda]_p$
and define $c_k$ for all $k$ by enforcing the rule
$c_{k+\lambda}=c_k$. 
Supposing that the desired conclusion does not hold, one has
$$p^{\lambda-i}[c_1,\dots,c_\lambda,\underbrace{p-1,\dots,p-1}_{i}]_p\equiv [c_1,\dots,c_\lambda,\underbrace{p-1,\dots,p-1}_{i},\underbrace{0,\dots,0}_{\lambda-i}]_p$$
$$
\equiv
[c_1,\dots,c_\lambda]_p+[\underbrace{p-1,\dots,p-1}_{i},\underbrace{0,\dots,0}_{\lambda-i}]_p
$$
$$
\equiv
[c_1,\dots,c_\lambda]_p-[\underbrace{0,\dots,0}_{i},\underbrace{p-1,\dots,p-1}_{\lambda-i}]_p$$
$$\equiv [c_{1+m},\dots,c_{\lambda+m}]_p=\bracket{p^mc}_q$$
for some integer $m$, where all the congruences are modulo $q-1$. It is impossible to have
$c_1=\cdots=c_i=0$  since this would force
the frequency of occurrence of the digit $0$
to differ in the digit strings $c_1,\dots,c_\lambda$
and $c_{1+m},\dots,c_{\lambda+m}$, which after all are just cyclic permutations one of the other.
Similarly we can rule out the possibility $c_{i+1}=\cdots=c_\lambda=p-1$.
Thus the base $p$ expansion of $c$ takes the form
$$c=[\underbrace{0,\dots,0}_{\alpha},
\underbrace{\bullet,\dots,\bullet}_{\beta},
\underbrace{p-1,\dots,p-1}_{\gamma}
]_p,$$
where 
$$\alpha< i,\;\;\beta>0,\;\;\;\gamma<\lambda-i,\;\;\;
\alpha+\beta+\gamma=\lambda,$$
and the bullets hold the place of a digit string not beginning with a $0$ and not ending with a $p-1$.
Then one has
$$\begin{array}{rcl}
1+\kappa_p(c)&=&(c+1)/p^\gamma\\
&>&(c+1-p^{\lambda-i})/p^\gamma+1\;\;(\mbox{strict inequality!})\\
&\geq &\tau_p(c+1-p^{\lambda-i})+1\\
&=&\tau_p(\langle p^mc\rangle_q)+1\\
&\geq &\kappa_p(\langle p^mc\rangle_q)+1
\end{array}
$$
in contradiction to the preceding lemma. This contradiction finishes the proof. \qed

 \end{document}